\documentclass{amsart}

\usepackage[all]{xy}

\newtheorem{theorem}{Theorem}[section]

\theoremstyle{definition}
\newtheorem{definition}[theorem]{Definition}
\newtheorem{example}[theorem]{Example}

\newtheorem{corollary}[theorem]{Corollary}

\theoremstyle{remark}
\newtheorem{remark}[theorem]{Remark}

\numberwithin{equation}{section}

\begin{document}
\author{Yanlong Hao}
\address {School of Mathematical Sciences,
         Tianjin 300071, P.R.China}
\email{haoyanlong13@mail.nankai.edu.cn}

\author{Xiugui Liu$^*$}
\address{School of Mathematical Sciences and LPMC,
         Nankai University,
         Tianjin 300071, P.R.China}
\email{xgliu@nankai.edu.cn}
\thanks{The second author was supported in part by the National Natural Science Foundation of China (No. 11171161), Program for New Century Excellent Talents in University (No. NCET-08-0288), and the Scientific Research Foundation for the Returned Overseas Chinese Scholars, State Education Ministry.}

\date{}
\title{On the rational homotopical nilpotency index of principal bundles}
\maketitle
\textbf{Abstract.} Let \rm{Aut}($p$) denote the space of all self-fibre homotopy equivalences of a principal $G$-bundle $p: E\rightarrow X$ of simply connected CW complexes with $E$ finite. When $G$ is a compact connected topological group, we show that there exists an inequality
$$n-{\rm N}(p)\leq {\rm Hnil}_{\mathbb{Q}}({\rm{Aut}}(p)_0)\leq n$$
for any space $X$, where $n$ is the number of non-trivial rational homotopy groups of $G$ and ${\rm N}(p)$ is defined in Section 2. In particular, ${\rm Hnil}_{\mathbb{Q}}({\rm{Aut}}(p)_{0})=n$ if $p$ is a fibre-homotopy trivial bundle and X is finite.

\section{Introduction}

Let $G$ be a compact connected topological group, and $p: E\rightarrow X$ be a principal  $G$-bundle of connected CW complexes with $E$ finite.
The identity component ${\rm{Aut}}(p)_0$ of the monoid of self-fibre homotopy equivalences of $E$ is known to be a grouplike space of CW type \cite[Prop. 2.2]{Y2}.

We first recall the general definition of homotopy nilpotency index for grouplike
spaces introduced by Berstein and Ganea \cite{I}.

Suppose $(G,\mu)$ is a grouplike
space with homotopy inverse $\nu: G\rightarrow G$. Let
\begin{equation*}
\varphi_2: G\times G\rightarrow G \end{equation*} denote the
commutator map
\begin{equation*}
\varphi_2(g_1,g_2)=\mu(\mu(g_1,g_2), \mu(\nu(g_1),\nu(g_2))).
\end{equation*}
Extend this definition to $\varphi_n: G^n\rightarrow G$ by
the rule
\begin{equation*}
\varphi_n=\varphi_2(\varphi_{n-1}\times id_G).
\end{equation*}
for $n>2$.

The \emph{homotopical nilpotency index} Hnil($G$) is the least integer (or infinite)
such that $\varphi_{n+1}$ is nullhomotopic. The \emph{rational homotopical nilpotency
index} ${\rm{Hnil}}_{\mathbb{Q}}(G)$ is set to be Hnil($G_{\mathbb{Q}}$) where $G_{\mathbb{Q}}$ denotes classical rationalization.

In \cite[Thm. 5.2]{Y2}, it is proven that ${\rm{Hnil}}_{\mathbb{Q}} ({\rm{Aut}}(p)_0)\leq n$, where $n$ is the
number of non-trivial rational homotopy groups of $G$. Here we find a lower
bound which in particular generalizes the main result of \cite{X}.
\begin{theorem}\label{5}
Let $p: E\rightarrow X$ be a principal
$G$-bundle of simply connected CW complexes with
classifying map $f: X\rightarrow B_{G}$,
where $E$ is finite. Then we have the following inequality:
\begin{equation*}
n-{\rm{N}}(p)\leq {\rm
Hnil}_{\mathbb{Q}}({\rm{Aut}}(p)_0)\leq n,
\end{equation*}
where $n$ is the number of non-trivial rational homotopy groups of $G$ and the quantity ${\rm{N}}(p)$ is defined in Section 2.
\end{theorem}
In particular, we have
\begin{corollary}\label{22}
When the fibre bundle is trivial, ${\rm Hnil}_{\mathbb{Q}}({\rm{Aut}}(p)_0)=n$.
\end{corollary}

\section{Background, tools and proofs}
For the sake of completeness, we first review some background and notations.
We assume that the reader is familiar with the basics of
rational homotopy theory. The book \cite{Y1} is a standard and excellent reference.

In \cite{D}, Sullivan defined a functor $A_{PL}(-)$ from
topological spaces to commutative differential graded algebras over
$\mathbb{Q}$ (CDGA for short) that is connected to the
cochain algebra functor $C^*(-;\mathbb{Q})$ by a sequence of natural
quasi-isomorphisms. Let $\wedge V$ denote the free commutative
graded algebra on the graded rational vector space $V$. A CDGA
$(A, d)$ is a {Sullivan algebra} if $A\cong \wedge V$ and $V$
admits a basis $(v_i)$ indexed by a well ordered set such that
$d(v_i)\in\wedge (v_j, j<i)$. If the differential $d$ is a
decomposable differential, that is, $d(V)\subseteq \wedge^{\geq
2}V$, we say $(A, d)$ is {minimal}. A CDGA $(A, d)$ is a
Sullivan model for $X$ if $(A, d)$ is a Sullivan algebra and there
is a quasi-isomorphism $(A, d)\rightarrow A_{PL}(X)$. If $(A, d)$ is
minimal then it is the {minimal Sullivan model} of $X$.

Let $p: E\rightarrow X$ be a principal $G$-bundle of simply connected CW
complexes, with $G$ a compact connected topological group. Let $(A,d)$ be a CDGA model of $X$ and recall that the minimal model of $G$ is necessarily of the form $(\Lambda V,d)$ with $V$ finite dimensional and concentrated in odd degrees. Then, see \cite[\S15.f]{Y1}, the bundle $p$ is modeled by an inclusion $(A,d)\hookrightarrow (A\otimes\Lambda V, D)$ for which $DV\subset A$. In other words, $(A\otimes\Lambda V, D)$ is a relative Sullivan algebra with base algebra $(A,d)$.

We say that a set $\{v_i\}_{i\in I}$ generates the relative Sullivan algebra $(A\otimes\Lambda V, D)$, if $v_i\notin A$, $Dv_i\in A$ for each $i$, $A$ and the linear space $W$ spanned by $\{v_i\}_{i\in I}$ generate multiplicatively $(A\otimes\Lambda V, D)$.  For each generator of $(A\otimes\Lambda V, D)$, say $v\in W$, the cohomology class $\alpha_{v}\in H^*(A,d)=H^*(X,\mathbb{Q})$ represented by $Dv$ is called a characteristic class of the bundle $p$. These
classes characterize the rational homotopy type of the bundle as they encode
the rational homotopy type of the classifying map $X\rightarrow BG$.
\begin{remark}Note that the definition of the characteristic class $\alpha_{v}$
depends on the choice of generators of $(A\otimes\Lambda V, D)$.
For example, let $p$ be an $S^1\times S^3$-fibration with minimal Sullivan model
$$\varphi:(\Lambda e,0)\hookrightarrow (\Lambda(e,x,y),D),$$
where deg$(e)=2$, deg$(x)=1$, deg$(y)=3$, $Dx=e$ and $Dy=e^2$. On the one hand, $\{x,y\}$ is a set of generators of $(A\otimes\Lambda V, D)$ with $\alpha_{x}\neq 0$ and $\alpha_{y}\neq 0$, on the other hand, $\{x':=x, y':=y-ex\}$ is also a set of generators of $(A\otimes\Lambda V, D)$ with $\alpha_{x'}\neq 0$ and $\alpha_{y'}=0$.
\end{remark}
Let $W$ always be the linear space spanned by a set of generators of $(A\otimes\Lambda V, D)$.
Then we can define a linear map of degree 1,
$$f:W\rightarrow H^*(A,d),\ \ \ \ f(v)=[Dv].$$
Let $k$ be the number of dimensions $m$ such that $f^m:W^m\rightarrow H^{m+1}(A,d)$ is injective.
\begin{definition}\label{000}
N($p$) is the minimum of the different values $k$ resulting from different choices of generators of $(A\otimes\Lambda V, D)$.
\end{definition}

Now we recall two fundamental facts. On the one hand, based on \cite[Thm. 10]{U1}, it is proven in \cite[Thm. 1]{U2} that a Lie model of BAut($p)_{0}$ is given
by the differential graded Lie algebra (DGL henceforth) formed by the
positive degree $A$-derivations of $(A\otimes \Lambda V,D)$. Recall that an $A$-derivation
$\theta\in {\rm{Der}}^p_{A}(A\otimes \Lambda V)$ of degree $p$ is a self linear map of degree $-p$ of $A\otimes \Lambda V$
which vanishes on $A$ and satisfies $\theta(xy) = \theta(x)y +(-1)^{p|x|}x\theta(y)$ for $x,y\in A\otimes \Lambda V$. The graded
vector space ${\rm{Der}}^p_{A}(A\otimes \Lambda V)$
 is a DGL with the commutator bracket
$[\theta_1, \theta_2] = \theta_1\theta_2-(-1)^{|\theta_1||\theta_2|}\theta_2\theta_1$ and differential $D(\theta) = [D,\theta]$. In particular,
$$\pi_*({\rm{Aut}}(p)_0)\otimes\mathbb{Q}\cong H_*({\rm{Der}}_{A}(A\otimes \Lambda V))$$
as Lie algebras, where we are considering the Samelson bracket in $\pi_*$(Aut$(p)_0)\otimes\mathbb{Q}$.
On the other hand, in \cite[Prop. 2.3]{Y2} it is proven that
$${\rm{Hnil}}_{\mathbb{Q}}({\rm{Aut}}(p)_0)={\rm{nil}}(\pi_*({\rm{Aut}}(p)_0))$$
the latter being the usual nilpotency index of a Lie algebra. In conclusion,
$${\rm{Hnil}}_{\mathbb{Q}}({\rm{Aut}}(p)_0)={\rm{nil}}(H_+({\rm{Der}}_{A}(A\otimes \Lambda V))).$$

Now we are in a position to prove Theorems \ref{5}.
\begin{proof}[Proof of Theorem \ref{5}]
In \cite[Theorem 5.2]{Y2}, it is proven that ${\rm{Hnil}}_{\mathbb{Q}}({\rm{Aut}}(p)_0)\leq n.$

Write $l=n-{\rm{N}}(p)$. By definition, there exists a set of generators of $A\otimes \Lambda V$ such that there are $l$ different dimensions, say $m_{1}<\cdots<m_{l}$, for which $f^{m_{i}}$ is not injective. In other words, for each $i=1,2,\ldots,l$ there is a non zero element $v_{i}\in W^{m_i}$ such that $Dv_{i}=d\Phi_i$ with $\Phi_i\in A$. Then, replace $v_{i}$ by $v_{i}-\Phi_i$ so that $Dv_{i}=0$ and we will assume that henceforth.
Enlarge $v_{1}$, $\ldots$, $v_{l}$ to a set of generators of $A\otimes \Lambda V$, denote $v_{0}=1$ and proceed as follows:
For each $0\leq i<j\leq l$, define $\theta_{j,i}\in{\rm{Der}}_A(A\otimes \Lambda V )$ as the $A$-derivation
which sends $v_j$ to $v_i$ and vanishes in any other element of the fixed basis of $V$.

It is immediate to see that $\theta_{j,i}$ is a cycle and not a boundary (as $v_j$ will
never be a boundary). Moreover, a straightforward computation shows that if $i<j<h$,
$$[\theta_{h,j},\theta_{j,i}]=\theta_{h,i}.$$
In particular,
$$[[\ldots[\theta_{l,l-1},\theta_{l-1.l-2}],\theta_{l-2,l-3}],\ldots],\theta_{1,0}]=\theta_{l,0}.$$
The desired result follows.
\end{proof}
\begin{proof}[Proof of Corollary \ref{22}]
When the fibre bundle is trivial, ${\rm{N}}(p)=0$. Hence we have $${\rm
Hnil}_{\mathbb{Q}}({\rm{Aut}}(p)_0)=n.$$
\end{proof}
\begin{example}
Let $p:E\rightarrow \mathbb{C}P^m$ be a principal $G$-bundle with $G$ a compact connected topological group. We show that N$(p)\leq 1$, which means
$$n-1\leq{\rm{Hnil}}_{\mathbb{Q}}({\rm{Aut}}(p)_0)\leq n.$$
If N$(p)\neq 0$, let $k$ be the smallest integer such that
$f:W^k\rightarrow H^{k+1}(\mathbb{C}P^m,\mathbb{Q})$
is injective. So $W^k$ is generated by one element $v$ and $f(v)=a^i$ where $a$ is a generator of $H^{2}(\mathbb{C}P^m,\mathbb{Q})$. For each generator of $W^{>k}$, say $w$, we have $f(w)=\lambda a^j$. Replace $w$ by $w-\lambda a^{j-i}v$ so that $f(w)=0$, so N($p)\leq1$, which means N($p)=1$.

In fact, it is easy to show that in this case ${\rm{Hnil}}_{\mathbb{Q}}({\rm{Aut}}(p)_0)=n-1$ if and only if N$(p)=1$ and the dimension $k$ used above is not the top dimension of the rational homotopy groups of $G$.
\end{example}
\begin{example}
When $G=S^1$, the model of a $G$-fibration is just $(A\otimes \Lambda v,D)$. Independently of the value of $N(p)$, we have
$$Hnil_Q(Aut(p)_0)=1,$$
since there is a non-trivial derivation that sends the generator $v$ to $1$.
\end{example}

\end{document}